# Dynamic Origin-Destination Estimation Using Smart Card Data: An Entropy Maximisation Approach

Abderrahman Ait Ali [a,1], Jonas Eliasson [b]
[a] The Swedish National Road and Transport Research Institute (VTI)
P.O. Box 55685, 114 28 Stockholm, Sweden
[1] E-mail: abderrahman.ait.ali@vti.se, Phone: +46 (0) 8 555 367 81
[b] Department of Science and Technology, Linköping University
Luntgatan 2, SE-602 47 Norrköping, Sweden

**Abstract**
Problems of dynamic origin-destination (OD) estimation using smart card data can be modelled using entropy maximisation and solved for large networks using solution techniques such as Lagrangian relaxation. In this paper, we give an overview of the research literature about OD estimation. We show how entropy maximisation can be used to model this problem in case station-entry data from smart cards is the only available information, i.e. number of entries is known but not exits nor the flow between the stations.
The large entropy maximisation program is solved using Lagrangian relaxation. The model is tested on a case study from the commuter train service in Stockholm with smart card entry-data from a working day in 2015. The results show that, given the entry-data, the entropy maximisation-based model and methods allows to find an OD-estimate that is as accurate as the reported estimates from manual observations. The estimates are further improved when using relevant assumptions and additional available data.

## 1  Introduction

The problem of origin-destination (hereafter OD) estimation appears in various domains and has been extensively studied by many researchers. Figure 1 shows the four-step-based transport planning model, OD estimation forms the second step, i.e. trip distribution. It is an important component and often comes after the trip generation where the zones (i.e. origins and destinations) are determined by their travel demand. The output of the OD estimation can be later used for mode choice and route assignment (McNally 2008).

The increasing availability of data (i.e. *big data*) led to the development of many new estimation models that can improve the quality of the trip distribution estimation in transport planning models. Many public transport networks generate valuable data through systems such as *automated fare collection* (AFC) using travel cards, referred to as *smart cards*. This data can be useful for transportation planner in different levels of planning, i.e. strategic, tactical or operational planning (Pelletier, Trépanier, and Morency 2011). Thus, a good combination of models and data would allow more accurate estimates of travel demand for a more optimised supply of transport infrastructure and operation.

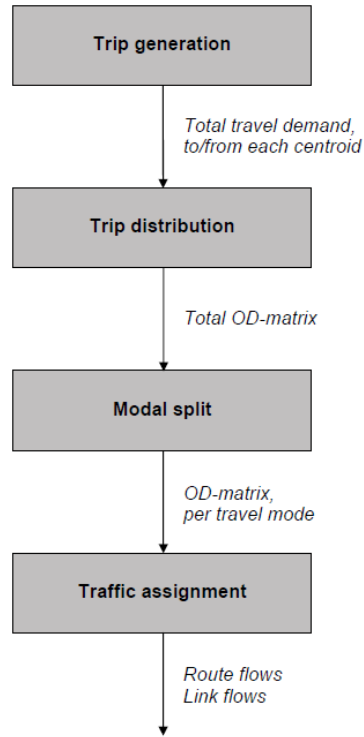

Figure 1: The four-step transport planning model, figure from (Peterson 2007)

Unlike *entry-exit* public transport systems, *entry-only* systems require the estimation of trip destinations from the collected smart card data. Travellers in entry-exit systems have to use their cards before entering their origin station and before exiting their destination station. The latter is not the case for entry-only transport systems where the card is only used before entering the system. Hence, data records about the trip destinations as well as the flow between the OD-pairs are not available and need to be estimated.

The aim of this paper is to develop and use entropy maximisation to estimate the trip destinations given entry-only smart card data. We review the research literature and start from a well-established entropy-maximisation formulation, we derive solution methods based on Lagrangian relaxation with symmetric trip assumption and additional available data. The case study is from the commuter train service in Stockholm which is an entry-only transport system, with smart card data from a working day in 2015.

The paper starts in section 2 with a brief literature review of the studied OD-estimation problems as well as the commonly used models. The mathematical formulation of the model and the solution methods are presented in section 3. The case study and the results are given in section 4. Section 5 concludes the paper.

## 2 Literature Review

In this section, we provide an overview of the research literature treating different OD estimation problems and models. In this review, we present the general formulation of the OD estimation problem and the adopted entropy maximisation approach.

### 2.1 OD Estimation Problem(s)

The research literature on OD estimation problems includes three main methodologies: data survey-based methods, trip distribution models and OD-matrix updating methods (Doblas and Benitez 2005). In the first, the estimation is based on a data survey of the trip subjects (e.g. survey of car-owners, public transport travellers, etc.). The second is based on estimation models that use incomplete or basic trips data (e.g. link flows in roads). The last uses various relevant information (e.g. growth rates) to update an old OD estimate. Depending on the needs and the available data in the case study, the methods and techniques to use can be a combination of one or more of these methodologies.

**Static and Dynamic**
There are two general variants of the OD estimation problem: *static* and *dynamic* (Deng and Cheng 2013). The static variant, also called *time-independent*, focuses on finding an estimate of the total number of trips between pairs of zones (e.g. stations) over a certain time interval (e.g. year, day or peak hours). The dynamic variant, also called *time-dependent*, adds the time dimension and therefore attempts to find an estimate of the number of trips in every time period (e.g. minutes or hours) over a certain time interval (e.g. a working day). The latter is more complex and requires larger sets of observation data and the former can be seen as a special case of the latter.

Many studies looked at the static OD estimation (Wang, Gentili, and Mirchandani 2012) whereas fewer others dealt with the dynamic version of the problem (Cho, Jou, and Lan 2009), some others have studied both (Deng and Cheng 2013). This paper will focus on the latter version, i.e. dynamic or time-dependent OD matrix estimation. The methods developed are valid even for the static problems.

**Zones and Traffic**
The OD estimation problems can also differ in terms of the considered zones and the studied type of transport traffic (i.e. agents). The type of studied traffic flow (cars, trains, train passengers, pedestrians, etc.) is determinant in formulating the OD estimation problems. Some studies such as (Wang, Gentili, and Mirchandani 2012) considered road vehicles as traffic agents and road sensor locations as the zones. This study focuses on train passenger flow in a commuter system and has therefore train stations as the trip zones and the passenger as traffic agents.

**Target Matrix**
Many formulations of the problem include the so-called *target matrix*, noted $\tilde{n}$ in equation (1). It is a reference matrix for comparison with the estimated one which can be seen as an updated version of the target matrix. Equations (2) and (5) illustrate how the target matrix is typically used in the growth OD estimation models and minimum information models, respectively.

In practice, old observations of the OD matrix are stored in order to be used as target matrices. These are used in OD estimation models attempting to find what the future OD matrices look like. The literature include studies that make use of target matrices (Wang and Zhang 2016) and other not using any (Cho, Jou, and Lan 2009). Like the latter studies, we do not use here a target matrix due to the absence of reliable data.

**General Formulation**
There are several variants of the OD estimation problem that can be covered by the general formulation presented in this section. The aim is generally to find an estimate of a matrix $n$, i.e. number of trips between OD pairs in the studied (transportation) system. This estimate induces a certain link flow $f$ that can be fully but often partially observed. The accuracy of the estimation depends, among other, on the observed available information describing the flows at certain links or zones.

Given certain observations of the link flows $\tilde{f}$ and a target OD-matrix $\tilde{n}$, the objective function in the general formulation is to minimise the deviation of the estimates from the observations. Equation (1) provides a general formulation where $d_n$ and $d_f$ are functions to measure the deviation of the OD-matrix and the link flows, respectively. The parameters $\alpha_n$ and $\alpha_f$ are weights allowing, for instance, to adjust the uncertainties on the deviations.

$$\min_{n,f} \alpha_n d_n(n, \tilde{n}) + \alpha_f d_f(f, \tilde{f}) \ . \tag{1}$$

Several studies have adopted this general formulation. For instance, Xie, Kockelman, and Waller (2011) used an entropy function for $d_n$ and adopted a least square estimator for $d_f$. In this study, we adopt the general formulation where we only consider the OD-matrix deviation function, i.e. link flows are ignored due to the absence of data. The smart card observations are used as constraints. Due to the inexistence of data on the target matrix $\tilde{n}$, we use entropy as a measure of the OD deviation $d_n$ as explained further in the paper.

## 2.2 OD Estimation Models

**Growth and Fratar Models**
Growth factor models refer to OD estimation methods where an outdated matrix is updated based on the estimation and calibration of different growth rates (Evans 1970). These methods are also simply called *growth models*. Equation (2) shows how an outdated matrix $\tilde{n}_{ij}$, of trips from $i$ to $j$, is updated using a *growth factor* $f$ to get an estimate of the new OD-matrix $n_{ij}$.

$$n_{ij} = f \, \tilde{n}_{ij} \ . \tag{2}$$

The growth factor $f$ can be a constant or an average rate, e.g. computed from survey data (Willekens 1983). It can also be compute by iteratively balancing factors $f_i$ and $f_j$ (origin-depend and destination-dependent growth factors, respectively) so that $f = f_i f_j$ where $f$ is defined as in equation (2). This growth variant is commonly called the *Furness method* (Morphet 1975).

**Gravity Models**
Gravity models attempt to estimate the OD-matrix using a gravitational attraction model, hence the name. Initially introduced as *Reilly's law* of retail gravitation where the main idea was that the attractiveness of retails depends on the physical distance and the size of the retail centres (Reilly 1931). This idea was used later to develop a mathematical model to predict traffic patterns (such as OD estimates) based on land-use (Voorhees 1955).

A possible formulation for gravity models is given in equation (3) where $c_{ij}$ is a *deterrence function* which relates to the generalized travel cost and $G_{ij}$ is a parameter that depends on the OD-pair, $f_i$ and $f_j$ are balancing factors as in the Furness method.

$$n_{ij} = G_{ij} \frac{f_i f_j}{c_{ij}} \qquad (3)$$

**Maximum Entropy and Minimum Information**
The maximum entropy approach is based on the statistical theory of probability and is proven to be equivalent to the gravity models (Wilson 1967). This approach relies on the idea that there are many possible trip distributions (i.e. *system states*) and that the most probable state (i.e. OD estimate) is the one that maximises the total entropy (randomness).

Wilson (1967) showed using the field of *statistical mechanics* that the total entropy $S_{entropy}$ of a system, subject to OD estimation, can be formulated as in equation (4) where log is the natural logarithm and $n = (n_{ij}^t)$ is the dynamic or time-dependent OD matrix estimate. This same formulation was also adopted in more recent studies such as (Xie, Kockelman, and Waller 2011).

$$S_{entropy}(n) = \sum_{ijt} (n_{ij}^t \log(n_{ij}^t) - n_{ij}^t) \qquad (4)$$

A similar (and equivalent) approach has been developed based on the *information theory*. Equation (5) presents a formulation of the deviation function in this approach. It uses the principle of *minimum information* instead of maximum entropy to find an OD-matrices based on a reference information on the trip distribution, i.e. OD target matrix. Van Zuylen and Willumsen (1980) were among the firsts to use this principle to build an OD estimation model. Instead of this approach, we use the entropy approach since we do not have reliable data for $\tilde{n}_{ij}^t$, i.e. time-dependent target matrix.

$$I(n, \tilde{n}) = \sum_{ijt} n_{ij}^t \log\left(\frac{n_{ij}^t}{\tilde{n}_{ij}^t}\right) \qquad (5)$$

**Other Models**
There are several studies that used various other models and additional techniques to both model and solve the OD estimation problem. For instance, logit and discrete choice models were used to estimate trip distributions by combining the destination choice and mode choice models (Ben-Akiva 1985). Different statistical techniques have been adopted such as generalized least square (Cascetta and Nguyen 1988), Bayesian inference (Maher 1983) and principal component analysis (Djukic et al. 2012). More recently, Kalman filters (Cho, Jou, and Lan 2009) and Markov chain models (Abareshi, Zaferanieh, and Safi 2019) were also used. These additional models show that the OD estimation problem has been studied from various perspectives using different tools hence the richness of its research literature.

# 3 Trip Destination Estimation

## 3.1 Basic Entropy Maximisation

**Motivation**

As mentioned in the literature review, entropy maximisation models were shown to be equivalent to gravity models (Wilson 1967) and that they are related to the minimum information principle as well as maximum utility models such as discrete choice models (Anas 1983). Entropy maximisation approach is therefore adopted and developed in this paper. Moreover, the choice of this approach can be motivated by the problem setup and the available input data as mentioned in the literature review.

**Methodology**

Entropy maximisation models can be formulated as convex optimization problems where the objective function to maximise is the total entropy of the system (as in equation 4) under certain (linear) constraints. Solving such as model for real world instance such as time-dependent and large networks is generally hard. Alternative solution methods to this optimization problem use Lagrangian relaxation which can give an estimate of the entropy maximising dynamic OD matrix. Once the optimization problem is formulated, the constraints are relaxed and associated to the Lagrangian multipliers. Using the KKT optimality conditions on the Lagrange function, we find an expression of the primal solution, i.e. OD matrix estimate. Additional data and assumptions can be included to get more accurate estimates.

**Notations**

Table 1 presents the notations that were adopted in the formulations that follow.

Table 1: Adopted notations

| Notation | Meaning |
|---|---|
| $T$ | set of time intervals over a working day |
| $S$ | set of all the nodes, i.e. train stations |
| $|S|$ | number of possible destinations for any given origin. |
| $n_{ij}^t$ | number of trips from station $i$ to $j$ at time interval $t$ |
| $O_i^t$ | number of passengers entering station $i$ at $t$ |

The smart card input data is provided by the values $O_i^t$ for all $i \in S$ and $t \in T$, i.e. the number of passengers *entering* an *origin* stations $i$ during a time interval $t$ and we want to find the values of the unknown variables $n = (n_{ij}^t)$. Hence, the distribution of the trip destination $D_j^t \coloneqq \sum_i n_{ij}^t$, number of passengers *exiting* a *destination* station $j$ during a time interval $t$.

To simplify the formulations, we write the indices $i$ without the inclusion $\in S$ and $t$ without $\in T$. Furthermore, since there is no trip from and to the same station, i.e. $n_{ij}^t$, we omit these trips in the summations, i.e. $\sum_{ijt} n_{ij}^t$ instead of $\sum_{ijt, i \neq j} n_{ij}^t$.

**Basic Model**

Using the previous notations, the formulation of the basic entropy maximisation model is given in equation (6). The set of constraints in equation (7) enforces that the sum of all trips originating from a station $i$ during a time interval $t$ should be equal to the smart card data $O_i^t$, i.e. collected number of entries to station $i$ during the time interval $t$.

$$\max_{n_{ij}^t \geq 0} \sum_{ijt} (n_{ij}^t \log(n_{ij}^t) - n_{ij}^t). \tag{6}$$

$$\sum_j n_{ij}^t = O_i^t \; ; \forall i, t. \tag{7}$$

The mathematical program is a nonlinear convex minimisation (concave maximisation) problem. For real world scenarios with many train stations, it is hard to solve using existing solvers. Hence, the use of Lagrangian relaxation to simplify the solution method. Thus, the sets of constraints in equation (7) are relaxed with associated Lagrangian multipliers ($\lambda_{it}$). This leads to the Lagrange function $L$ in equation (8), i.e. relaxed dual objective function.

$$L(n, \lambda) = \sum_{ijt} (n_{ij}^t \log(n_{ij}^t) - n_{ij}^t) + \sum_{it} \lambda_{it} (\sum_j n_{ij}^t - O_i^t). \tag{8}$$

For a fixed value of $\lambda$, $L$ can be minimised with respect to the variable $n$ and the optimality conditions provides an expression of the OD estimate. Equation (9) presents the primal solution as a function of the dual multipliers. Enforcing the set of constraints (7) allows to find an expression for the multipliers $\lambda$. Equation (10) presents this step as well as the formulation of the entropy maximising OD estimate as a function of the smart card entry counts.

$$\frac{\partial L}{\partial n_{ij}^t} = 0 \Rightarrow n_{ij}^t = e^{\lambda_{it}}. \tag{9}$$

$$\sum_j n_{ij}^t = O_i^t \Rightarrow \sum_j e^{\lambda_{it}} = O_i^t \Rightarrow n_{ij}^t = \frac{O_i^t}{|S|}. \tag{10}$$

This formulation is what we refer to in the case study as the basic method (BM). It is basic since the trips are distributed uniformly among all the destination stations. This means that the probabilities of traveling to a certain destination station is the same (i.e. $\frac{1}{|S|}$) regardless of the destination station. This is not likely to happen in reality because different stations have different attractivities at different times and these probabilities should therefore be different from one destination to another.

### 3.2 Symmetry Assumption

To reflect that the destination station may have an influence on the trip destination probability, we assume that the total trips over the day are *symmetric*. This means that all trips originating from an origin station for work (or leisure) are also repeated in the opposite direction (to go back home to the same origin station some time during the day). Equation (11) formulate this symmetric assumption as a set of additional constraints to the entropy maximisation model.

$$\sum_{it} n_{ij}^t = \sum_t O_j^t \ ; \ \forall j. \tag{11}$$

Relaxing these constraints, with $\mu_j$ as the corresponding associated multipliers, and following the same Lagrangian relaxation steps, we find a new OD estimate. Equation (12) gives an expression this new estimate.

$$\frac{\partial L}{\partial n_{ij}^t} = 0 \Rightarrow n_{ij}^t = e^{\lambda_{it}} e^{\mu_j}. \tag{12}$$

Using the set of constraints (7) and (11), we find the values of the multipliers which allow to satisfy both sets. Equation (13) shows this as well as the expression of the entropy maximising OD estimate given the symmetric assumption.

$$\begin{cases} \sum_j n_{ij}^t = O_i^t \\ \sum_{it} n_{ij}^t = \sum_t O_j^t \end{cases} \Rightarrow \begin{cases} e^{\lambda_{it}} \sum_j e^{\mu_j} = O_i^t \\ e^{\mu_j} \sum_{it} e^{\lambda_{it}} = \sum_t O_j^t \end{cases} \Rightarrow n_{ij}^t = O_i^t \frac{\sum_t O_j^t}{O}. \tag{13}$$

Unlike the estimate from the basic method, the symmetric assumption yields a distribution of the destination probabilities that differs from one destination station to another. These probabilities depend on the smart card counts of the total number of trips that are originating from the station. The attractivity of a destination station increases when more trips originate from it. This estimate is referred to as the symmetry assumption (SA) method. Note that if the total origin entries over the day are similar for all origin stations (i.e. $\sum_t O_j^t = \sum_t O_s^t$, $\forall j \neq s$), the OD estimate in equation (13) is reduced to the previous BM estimate in (10).

### 3.3 Additional Data

Given additional relevant data, it is possible to extend the OD estimation model and improve the accuracy of the estimated OD matrix. This data can come from various sources such as manual counting, passenger survey results or various travel prognosis providing additional information on, for instance, total average travelled distance or time, total person-km, ticket revenue or other financial results.

For the sake of illustration, we use the total reported person-km that we note $\bar{d}$ (e.g. from survey results). Equation (14) presents the additional constraint enforcing the total reported person-km where $d_{ij}$ is the travel distance between station $i$ and $j$.

$$\sum_{ijt} d_{ij} n_{ij}^t = \bar{d}. \tag{14}$$

The extended model includes additionally the constraint (14) and we note $\theta$ as the corresponding multiplier. Following the same steps and using the optimality conditions, we get the formulation of the new OD estimate in equation (15) where $A_{it} := e^{\lambda_{it}}$ and $e^{\mu_j}$ relates to the set of constraints (7) and (11), respectively. Determining the parameter $A_{it}$ using the set of constraints (7) leads to the formulation in equation (16) which presents the OD estimate using the additional data, referred to as AD method in the case study.

$$\frac{\partial L}{\partial n_{ij}^t} = 0 \Rightarrow n_{ij}^t = e^{\lambda_{it}+\mu_j+\theta k_{ij}} = A_{it} e^{\theta d_{ij}+\mu_j}. \tag{15}$$

$$\sum_j n_{ij}^t = O_i^t \Rightarrow A_{it} \sum_j e^{\mu_j} e^{\theta d_{ij}} = O_i^t \Rightarrow n_{ij}^t = O_i^t \frac{e^{\theta d_{ij}+\mu_j}}{\sum_{s \neq i} e^{\theta d_{is}+\mu_s}}. \tag{16}$$

Note that if the additional data is not considered (e.g. $d_{ij} = 0$), equation (16) is reduced using the set of constraints (11) to the symmetric assumption (SA) method.

The term $p(j|i) := \frac{e^{\theta d_{ij}+\mu_j}}{\sum_{s \neq i} e^{\theta d_{is}+\mu_s}}$ can be interpreted as the probability of choosing destination $j$ given origin $i$. In this case, the exponent $\theta d_{ij}$ is similar to a certain utility value ($\theta < 0$) for travelling between $i$ and $j$. This utility can include $m$ types of additionally relevant data $k_{ij}^{(1)}, \ldots, k_{ij}^{(m)}$ (if available) and form the total utility and choice probabilities as in equation (17).

$$u_{ij} = K_j + \theta_1 k_{ij}^{(1)} + \cdots + \theta_m k_{ij}^{(m)} \Rightarrow p(j|i) = \frac{e^{u_{ij}}}{\sum_{s \neq i} e^{u_{is}}} \Rightarrow n_{ij}^t = O_i^t p(j|i) \tag{17}$$

**Calibration**

The formulation in equation (17) is similar (no random error term) to a discrete choice model where the discrete choices here are the destination stations $j$ for each origin $i$. The parameters $K_j$ and $\theta_1, \ldots, \theta_m$ are specific to the (transport) system where the OD estimation is performed and need to be *calibrated* to reflect the (dis-)utilities of the travellers in the system. For instance, it is possible to estimate the values of these parameters using the additional data (or constraints), stated (or revealed) preference survey data or old OD matrices (target matrix) of the studied transport system.

In this study, we use the total person-km and calibrate the parameters $\theta$ and $\mu_j$ using the reported total person-km $\bar{d}$ expressed in constraint (14) and the symmetric assumption constraints in (11). The parameters are iteratively estimated until the satisfaction of the two constraints (up to a certain tolerance $\epsilon = 10^{-5}$). Figure 2 presents the iterative calibration algorithm, i.e. estimation of the parameters $\theta$ and $\mu_j$ in equation (16).

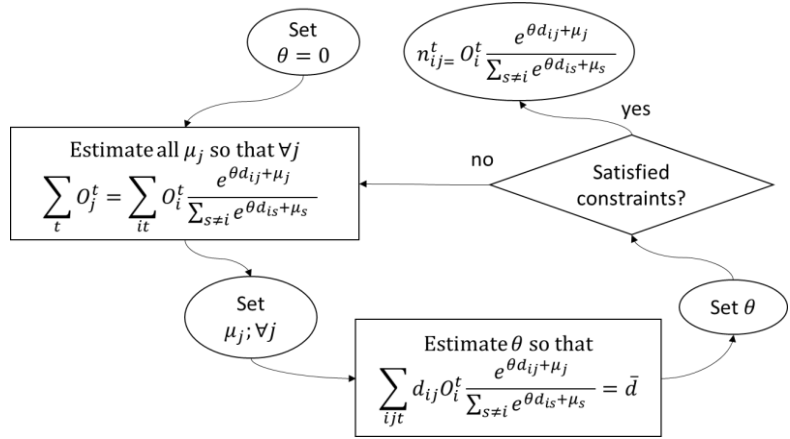

Figure 2: Iterative calibration algorithm for estimating the parameters $\theta$ and $\mu_j; \forall j$.

## 4  Case Study

This section describes how the model was implemented and tested. The data is described before presenting the results. Finally, the validation of the results is discussed.

### 4.1  Data

In this paper, we study the commuter train system in Stockholm, Sweden (locally called *pendeltåg*). Figure 3 illustrates the almost 250 km long commuter network of 2015 with around 58 stations. Due to the absence of accurate data for 7 minor stations , mostly on the short line between to/from *Gnesta*, these stations were neglected, i.e. $|S| + 1 = 51$ train stations.

The smart card data is from the locally called *SL Access Card* which is a travel smart card used by public transport travellers in Stockholm. Since the commuter train system is entry-only, we only have data on the number of passengers entering each station during different time periods of an entire working day in 2015 with 15 minutes time intervals, i.e. $|T| = 96$ intervals in total (for one working day). Figure 3 presents, in addition to the network, the total number of travellers entering each (origin) station. The neglected *Gnesta* line is also presented. Note that Stockholm Central station has the largest number of entries compared to the other stations, e.g. more than four times the second most frequented station.

Figure 4 illustrates the temporal distribution of the data that was used from smart cards. It shows the total number of passengers entering the commuter system in every 15min during the studied working day. Note the typical two peaks during the working days, one in the morning and a smoother one in the afternoon.

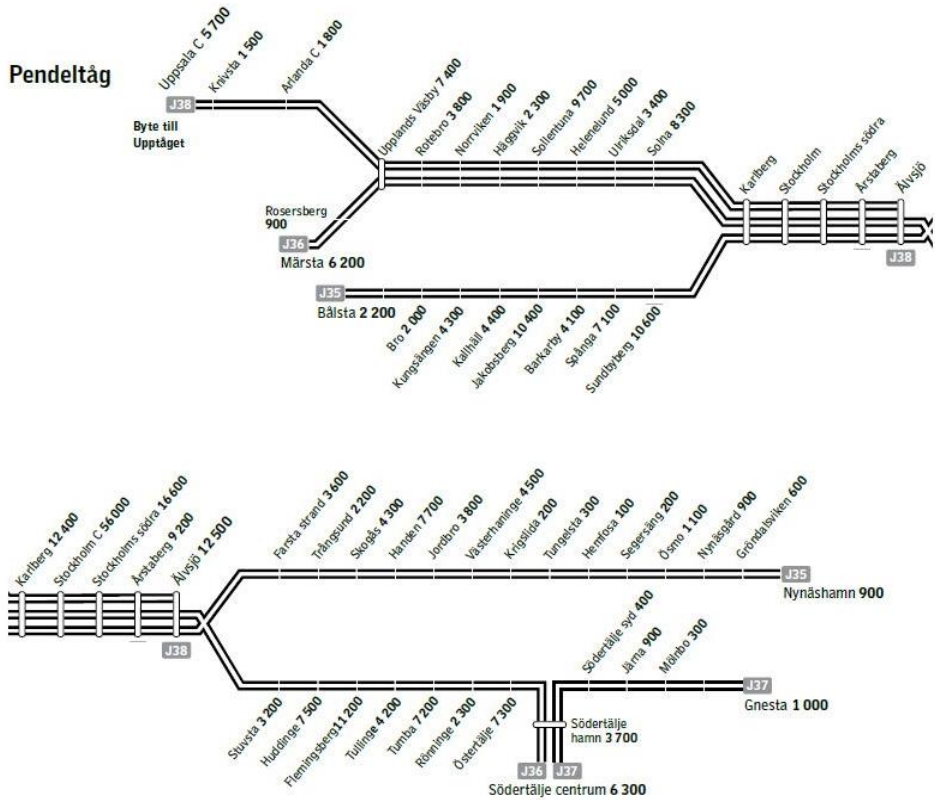

Figure 3: Stockholm commuter network in 2015 with total daily entries ($\sum_t O_i^t ; \forall i$) of passengers from the different origin stations (SLL, 2015)

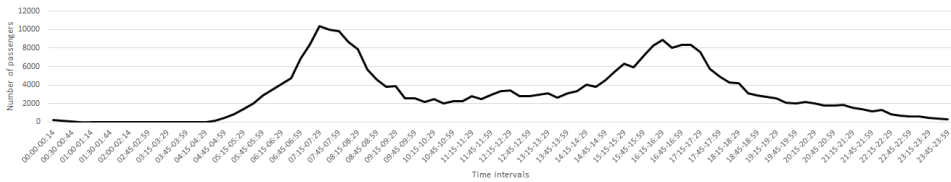

Figure 4: Total number ($\sum_i O_i^t ; \forall t$) of passengers entering to the commuter train system in every 15min time interval in a working day in Stockholm (SLL, 2015)

## 4.2 Tests and Results

We test the different OD estimation methods that were developed, i.e. basic method (BM), symmetric assumption (SA) and additional data (AD). This is done by estimating the dynamic OD matrix using the smart card data that was previously presented about Stockholm's commuter train system. Table 1 presents the OD methods that have been tested and their characteristics.

Table 1: OD estimation variants and their characteristics.

| Variant | Characteristics |
|---|---|
| BM | Basic method, only smart card data. |
| SA | BM with the symmetric assumption. |
| AD | SA with the additional data: distances and total person-km. |

We test and compare the estimation results from the different estimation variants. BM refers to the basic method in section 3.1, SA represents the method variant with the symmetric assumption as described in section 3.2 and finally AD is the extended method using additional available data for distances and total person-km as presented in section 3.3.

The various resulting OD estimates are computed and compared given the available data. Figure 5 compares the three (BM, SA and AD) estimated temporal variation of the distribution of the trip destination at Stockholm Central station, i.e. total number of exits from the station. The fourth (dashed) curve correspond to the smart card entry distribution, i.e. total number of entries to the station.

We retrieve the two peaks from the entry distribution in both SA and AD. However, the estimate from BM is flatter and gives (as expected) a worse estimate than the other methods.

The estimated number of passengers exiting at the central station is higher than the counts of entries to the station in the morning peak hours. In the afternoon peak hours and later, the opposite happens, i.e. higher counted entries than estimated exits. This remark is valid for all the tested methods (except BM). One can conclude from the estimation that the central station at Stockholm is used more as a departure station (higher entries) in the afternoon peak hours and more as a destination station in the morning peak hours.

AD yields generally higher estimates than SA (and BM). This can be explained by the fact that destination probabilities (or attractivities) are amplified using the additional data, i.e. distance compared to the other methods. This means that destination stations near the origin have higher attractivities which is the case for the central station.

A more general and detailed visualization of the resulting dynamic OD estimate can be made using a 3D surface plot. Figure 6 shows the total daily number of passengers (vertical axis) travelling from entry origin stations (right horizontal axis) to exit destinations (left horizontal axis). The peaks can be observed from/to Stockholm Central station due to the large number of passengers boarding from this station (according to smart card entry data). The largest passenger flow is between Stockholm Central (most frequented) and Southern station (second most frequented) which is noticeable while commuting in working days.

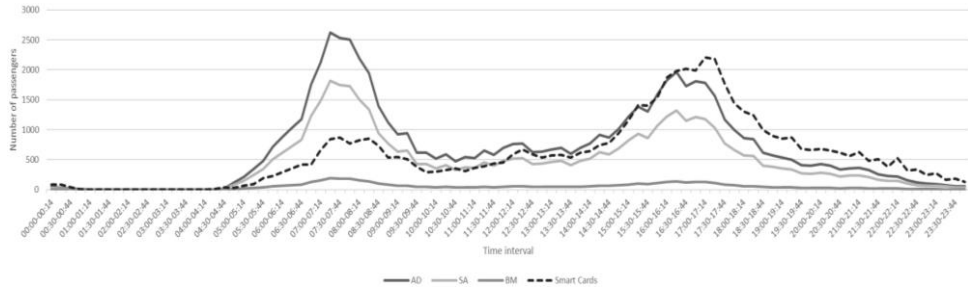

Figure 5: Comparison between the different trip destination estimates ($\sum_i n_{is}^t$; $s =$ central station, $\forall t$) as well as the trip entry smart card data ($O_s^t$; $\forall t$) at Stockholm Central station.

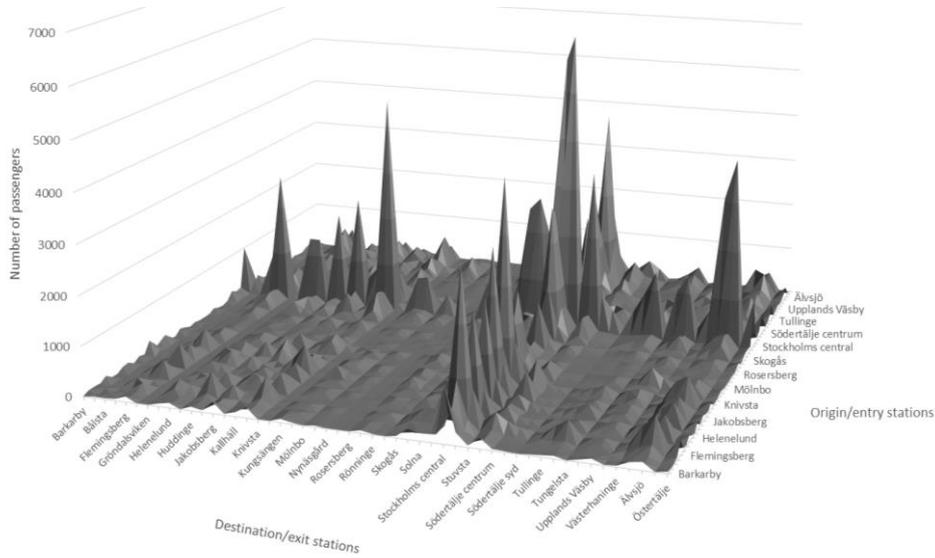

Figure 6: Resulting daily OD estimate ($\sum_s n_{ij}^t , \forall ij$) using additional data (AD).

### 4.3 Validation Discussions

In addition to the visualization of the OD estimates (figures 5 and 6), the methods can be validated by computing several aggregated statistics and compare them with available data from the local transport agency. Manual counting on platforms or stated preference methods are examples of methods that can be used to find these aggregated values.

Using the different OD estimates, we compute and compare travel statistics with the ones reported in 2015 by the transport agency *SLL* in Stockholm (SLL 2015). Table 2 presents the comparison between some of these travel statistics. The columns represent the different travel statistics that were compared and the rows show the estimated statistics for each of the studied OD-estimation variants, the last row shows the reported statistics from the local transport agency SLL.

Table 2: Comparison of OD-estimate travel statistics and the reported ones (SLL 2015).

| Variants | Total person-km (in $10^3$) | Average travel distance (in km) | Total daily exits at Central station |
|---|---|---|---|
| BM | 10 194 | 34.8 | 47 380 |
| SA | 7 178 | 24.5 | 45 293 |
| AD | 5 776 | 19.7 | 66 552 |
| SLL | 5 776 | 18.7 | 64 700 |

The travel statistics from the OD-estimate based on AD method are the closest to the ones reported by SLL (local agency). A part from the total person-km that was used as additional data during the estimation in AD, the average travel distance and the daily exits at the central station correspond to the reported values with an accuracy of 95% and 97%, respectively.

SA and, to a lesser extent, BM can yield results with an accuracy of up to 85% compared to the reported statistics. Therefore, an OD-estimate based on entropy maximisation may be used for OD estimation and can be as accurate as the currently used manual methods.

There are alternative ways to further validate the estimates. For instance, the use of automatic data such as door sensors to count exits and entries to the trains or the weight of trains between stations. Due to the lack of such data, we were not able to further compare and validate the results with these automatic sources of data.

## 5  Conclusions and Future Works

We have shown in the literature review that there are several methods for OD-matrix estimation. Entropy maximisation is the one that was used and developed in this study for modelling the dynamic OD estimation problem. We derived three different variants of the solution method based on Lagrangian relaxation: basic method (BA), symmetric assumption (SA) and additional data (AD). These were tested using real smart card data from the commuter train service in Stockholm in 2015. The results showed that the model and the solution methods, more particularly AD, allow to find an estimate of the real dynamic OD which yields results as accurate as the manual observations reported by the local transport agency.

The developed entropy maximisation methods showed that it is possible to get accurate OD-estimates based on automatic (cheap) sources of data instead of manual (expensive) ones. Additional available and inexpensive sources of data can be used to improve the estimation or further validation the resulting OD estimates.